# The Consistency, the Composition and the Causality of the Asynchronous Flows

Serban E. Vlad[1]
[1]Oradea City Hall, p-ta Unirii, nr. 1, 410100, Oradea, Romania.

**Abstract**

Let $\Phi:\{0,1\}^n \to \{0,1\}^n$. The asynchronous flows are (discrete time and real time) functions that result by iterating the coordinates $\Phi_i, i \in \{1,...,n\}$ independently on each other. The purpose of the paper is that of showing that the asynchronous flows fulfill the properties of consistency, composition and causality that define the dynamical systems. The origin of the problem consists in modelling the asynchronous circuits from the digital electrical engineering.

## 1. Introduction

The Boolean autonomous deterministic regular asynchronous systems have been defined by the author in 2007 and a study of such systems can be found in [12]. The concept has its origin in switching theory, the theory of modelling the asynchronous (or switching) circuits from the digital electrical engineering. The attribute Boolean vaguely refers to the Boole algebra with two elements; autonomous means that there is no input; determinism means the existence of a unique state function; and regular indicates the existence of a function $\Phi:\{0,1\}^n \to \{0,1\}^n$, $\Phi=(\Phi_1,...,\Phi_n)$ that 'generates' the system. Time is discrete: $\{-1,0,1,...\}$, or continuous: $\mathbf{R}$. The system, which is analogue to the (real, usual) dynamical systems, iterates (asynchronously) on each coordinate $i \in \{1,...,n\}$ one of

- $\Phi_i$: we say that $\Phi$ is computed, at that time instant, on that coordinate;

- $\{0,1\}^n \ni (\mu_1,...,\mu_i,...,\mu_n) \mapsto \mu_i \in \{0,1\}$: we use to say that $\Phi$ is not computed, at that time instant, on that coordinate.

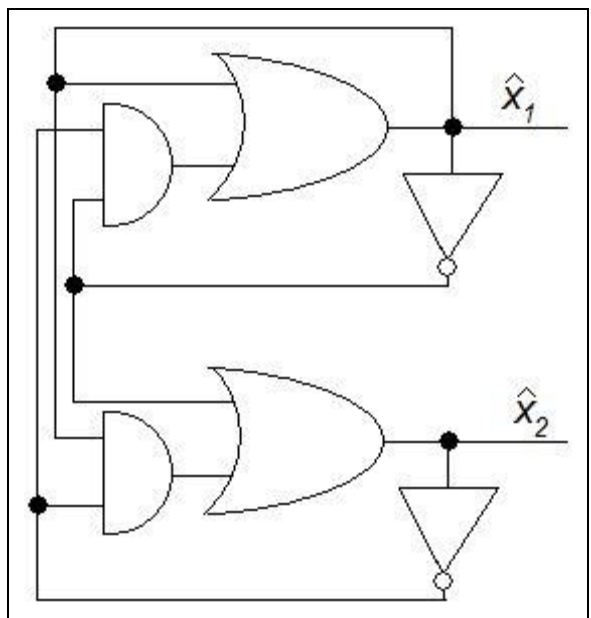

**Figure 1: Asynchronous circuit**

The flows are these that result by analogy with the dynamical systems.
The 'nice' discrete time and real time functions that the (Boolean) asynchronous systems work with are called signals. The functions that show when and how the coordinates $\Phi_i$ are computed are called computation functions. In order to point out the source of inspiration, we give the example of the circuit from Figure 1, where $\hat{x}:\{-1,0,1,...\} \to \{0,1\}^2$ is the signal representing the state of the system, and the initial state is $(0,0)$. The function that generates the system is $\Phi:\{0,1\}^2 \to \{0,1\}^2$, $\forall \mu \in \{0,1\}^2$,

$$\Phi(\mu) = (\mu_1 \cup \overline{\mu_1} \cdot \overline{\mu_2}, \overline{\mu_1} \cup \mu_1 \cdot \overline{\mu_2}).$$

The evolution of the system is shown in its state diagram from Figure 2, where the arrows indicate the time increase and we have underlined the coordinates $\mu_i, i=\overline{1,2}$ that, by the computation of $\Phi$, change their value: $\Phi_i(\mu) = \overline{\mu_i}$. Let $\alpha:\{0,1,2,...\} \to \{0,1\}^2$ be the computation function whose values $\alpha_i^k$ show that $\Phi_i$ is computed at the time instant $k$ if $\alpha_i^k = 1$, respectively that it is not computed at the time instant $k$ if $\alpha_i^k = 0$, where $i=\overline{1,2}$ and $k \in \{0,1,2,...\}$. The uncertainty related with the modelled circuit, depending in general on the technology, the temperature, etc, manifests in the fact that the order and the time of computation of each coordinate function $\Phi_i$ are not known.

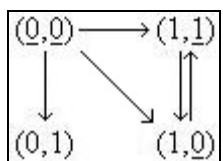

**Figure 2: The state diagram of the circuit from Figure 1**

The situation $\alpha^0 = (0,0)$, when no coordinate of $\Phi$ is computed at the time instant $0$, shows that the system remains in $(0,0)$.

If the first coordinate of $\Phi$ is computed at the time instant $0$, i.e. $\alpha^0 = (1,0)$, then Figure 2 indicates the transfer from $(0,0)$ in $(1,0)$.

We suppose that the second coordinate is computed at the time instant $0$, i.e. $\alpha^0 = (0,1)$, and in this case the system transfers from $(0,0)$ to $(0,1)$, where it remains indefinitely long for any values of $\alpha^1, \alpha^2, \alpha^3, ...$, since $\Phi(0,1) = (0,1)$. Such a signal $\hat{x}$ is called eventually constant and it corresponds to a stable system.

The last possibility is given by $\alpha^0 = (1,1)$ that indicates the transfer from $(0,0)$ to $(1,1)$, as resulted by the simultaneous computation of $\Phi_1(0,0)$ and $\Phi_2(0,0)$.

If the system is in one of the points $(1,0), (1,1)$ and the set $\{k \mid k \in \mathbf{N}, \alpha_2^k = 1\}$ is infinite, then it switches infinitely many times between $(1,0)$ and $(1,1)$ and this corresponds to an unstable system.

The purpose of our paper is that of showing that the flows of these systems fulfill the properties of consistency, composition and causality that define the dynamical systems.

## 2. Preliminaries. Signals

**Notation 1** We denote by $\mathbf{B} = \{0,1\}$ the binary Boole algebra. Its laws are the usual ones:

**Table 1.**

| $\overline{\phantom{x}}$ |     | $\cdot$ | 0 | 1 | $\cup$ | 0 | 1 | $\oplus$ | 0 | 1 |
|---|---|---|---|---|---|---|---|---|---|---|
| 0 | 1, | 0 | 0 | 0, | 0 | 0 | 1, | 0 | 0 | 1 |
| 1 | 0 | 1 | 0 | 1 | 1 | 1 | 1 | 1 | 1 | 0 |

and they induce laws that are denoted with the same symbols on $\mathbf{B}^n, n \geq 1$.

**Definition 2** Both sets $\mathbf{B}$ and $\mathbf{B}^n$ are organized as topological spaces by the discrete topology.

**Notation 3** $\mathbf{N}_- = \{-1, 0, 1, ...\}$ is the notation of the discrete time set.

**Notation 4** We denote
$$\hat{Seq} = \{(k_j) \mid k_j \in \mathbf{N}_-, j \in \mathbf{N}_- \text{ and } k_{-1} < k_0 < k_1 < ...\},$$
$$Seq = \{(t_k) \mid t_k \in \mathbf{R}, k \in \mathbf{N} \text{ and } t_0 < t_1 < t_2 < ... \text{ unbounded from above}\}.$$

**Notation 5** $\chi_A : \mathbf{R} \to \mathbf{B}$ is the notation of the characteristic function of the set $A \subset \mathbf{R}$: $\forall t \in \mathbf{R}$,
$$\chi_A(t) = \begin{cases} 1, \text{ if } t \in A, \\ 0, \text{ otherwise} \end{cases}.$$

**Definition 6** The *discrete time signals* are by definition the functions $\hat{x} : \mathbf{N}_- \to \mathbf{B}^n$. Their set is denoted with $\hat{S}^{(n)}$.

The *continuous time signals* are the functions $x : \mathbf{R} \to \mathbf{B}^n$ of the form $\forall t \in \mathbf{R}$,
$$x(t) = \mu \cdot \chi_{(-\infty, t_0)}(t) \oplus x(t_0) \cdot \chi_{[t_0, t_1)}(t) \oplus ... \oplus x(t_k) \cdot \chi_{[t_k, t_{k+1})}(t) \oplus ... \quad (1)$$

where $\mu \in \mathbf{B}^n$ and $(t_k) \in Seq$. Their set is denoted by $S^{(n)}$.

**Remark 7** The signals model the electrical signals of the circuits from the digital electrical engineering.

**Remark 8** At Notation 4 and Definition 6 a convention of notation has occurred, namely a hat ^ is used to show that we have discrete time. The hat will make the difference between, for example, the notation of the discrete time signals $\hat{x}, \hat{y}, ...$ and the notation of the real time signals $x, y, ...$

**Lemma 9** For any $x \in S^{(n)}$ and any $t \in \mathbf{R}, x(t-0) \in \mathbf{B}^n$ exists with the property
$$\exists \varepsilon > 0, \forall \xi \in (t - \varepsilon, t), x(\xi) = x(t - 0). \quad (2)$$

**Proof.** We presume that $x, t$ are arbitrary and fixed and that $x$ is of the form (1), with $\mu \in \mathbf{B}^n$ and $(t_k) \in Seq$. If $t \leq t_0$, then any $\varepsilon > 0$ makes (2) be true with $x(t-0) = \mu$; and if $k \geq 0$ exists with $t \in (t_k, t_{k+1}]$, then any $\varepsilon \in (0, t_{k+1} - t_k)$ makes (2) be true with $x(t-0) = x(t_k)$.

**Definition 10** The function $\mathbf{R} \ni t \mapsto x(t-0) \in \mathbf{B}^n$ is called the *left limit function* of $x$.

**Definition 11** The *discrete time forgetful function* $\hat{\sigma}^{k'} : \hat{S}^{(n)} \to \hat{S}^{(n)}$ is defined for any $k' \in \mathbf{N}$ by

$$\forall \hat{x} \in \hat{S}^{(n)}, \forall k \in \mathbf{N}_-, \hat{\sigma}^{k'}(\hat{x})(k) = \hat{x}(k+k') \tag{3}$$

and the *real time forgetful function* $\sigma^{t'} : S^{(n)} \to S^{(n)}$ is defined for $t' \in \mathbf{R}$ in the following manner

$$\forall x \in S^{(n)}, \forall t \in \mathbf{R}, \sigma^{t'}(x)(t) = \begin{cases} x(t), t \geq t', \\ x(t'-0), t < t'. \end{cases} \tag{4}$$

## 3. Computation functions

**Definition 12** The *discrete time computation functions* are by definition the sequences $\alpha : \mathbf{N} \to \mathbf{B}^n$. Their set is denoted by $\hat{\Pi}'_n$. In general, we write $\alpha^k$ instead of $\alpha(k), k \in \mathbf{N}$.

The *real time computation functions* $\rho : \mathbf{R} \to \mathbf{B}^n$ are by definition the functions of the form

$$\rho(t) = \alpha^0 \cdot \chi_{\{t_0\}}(t) \oplus \alpha^1 \cdot \chi_{\{t_1\}}(t) \oplus ... \oplus \alpha^k \cdot \chi_{\{t_k\}}(t) \oplus ... \tag{5}$$

where $(t_k) \in Seq$. Their set is denoted by $\Pi'_n$.

**Remark 13** The meaning of the computation functions $\alpha \in \hat{\Pi}'_n, \rho \in \Pi'_n$, subject to the additional property of progressiveness that will be stated later, is that of showing when –in discrete time and in real time- and how the Boolean functions $\Phi : \mathbf{B}^n \to \mathbf{B}^n$ are computed.

**Lemma 14** For any $\rho \in \Pi'_n$ and any $t \in \mathbf{R}$, we have

$$\exists \varepsilon > 0, \forall \xi \in (t-\varepsilon, t), \rho(\xi) = (0,...,0). \tag{6}$$

**Proof.** Analogue with the proof of Lemma 9.

**Definition 15** The *discrete time* $\hat{\sigma}^{k'} : \hat{\Pi}'_n \to \hat{\Pi}'_n, k' \in \mathbf{N}$ and the *continuous time* $\sigma^{t'} : \Pi'_n \to \Pi'_n$ *forgetful function*, $t' \in \mathbf{R}$, are defined by: $\forall \alpha \in \hat{\Pi}'_n, \forall k \in \mathbf{N}$,

$$(\hat{\sigma}^{k'}(\alpha))^k = \alpha^{k+k'} \tag{7}$$

and $\forall \rho \in \Pi'_n, \forall t \in \mathbf{R}$,

$$\sigma^{t'}(\rho)(t) = \rho(t) \cdot \chi_{[t',\infty)}(t). \tag{8}$$

**Remark 16** Definition 15, equation (8) was given by analogy with Definition 11, equation (4), taking into account (6): $\forall t \in \mathbf{R}$,

$$\sigma^{t'}(\rho)(t) = \begin{cases} \rho(t), t \geq t', \\ \rho(t'-0), t < t' \end{cases} = \begin{cases} \rho(t), t \geq t', \\ (0,...,0), t < t' \end{cases} = \rho(t) \cdot \chi_{[t',\infty)}(t)$$

indeed.

## 4. Progressiveness

**Definition 17** The discrete time computation function $\alpha \in \hat{\Pi}'_n$ is called *progressive* if

$$\forall i \in \{1,...,n\}, \text{the set } \{k \mid k \in \mathbf{N}, \alpha_i^k = 1\} \text{ is infinite}. \tag{9}$$

The set of the discrete time progressive computation functions is denoted by $\hat{\Pi}_n$.

The real time computation function $\rho \in \Pi'_n$ is called *progressive* if

$$\forall i \in \{1,...,n\}, \text{ the set } \{t \mid t \in \mathbf{R}, \rho_i(t) = 1\} \\ \text{ is unbounded from above} \tag{10}$$

is true. The set of the real time progressive computation functions is denoted by $\Pi_n$.

**Theorem 18** a) Let the computation function $\alpha \in \hat{\Pi}_n'$. The following equivalence holds:
$$\alpha \in \hat{\Pi}_n \Leftrightarrow \hat{\sigma}^1(\alpha) \in \hat{\Pi}_n.$$

b) The computation function $\rho \in \Pi_n'$ and $t' \in \mathbf{R}$ are given. The following equivalence holds:
$$\rho \in \Pi_n \Leftrightarrow \sigma^{t'}(\rho) \in \Pi_n.$$

**Proof.** a) For any $i \in \{1,...,n\}$, the sets $\{k \mid k \geq 0, \alpha_i^k = 1\}$, $\{k \mid k \geq 1, \alpha_i^k = 1\}$ are simultaneously finite or infinite.

b) We suppose that $\rho$ is of the form
$$\rho(t) = \rho(t_0) \cdot \chi_{\{t_0\}}(t) \oplus \rho(t_1) \cdot \chi_{\{t_1\}}(t) \oplus ... \oplus \rho(t_k) \cdot \chi_{\{t_k\}}(t) \oplus ... \quad (11)$$

with $(t_k) \in Seq$. We denote with $k' \geq 0$ the rank of the sequence $(t_k)$ that is defined by
$$\sigma^{t'}(\rho)(t) = \rho(t_{k'}) \cdot \chi_{\{t_{k'}\}}(t) \oplus \rho(t_{k'+1}) \cdot \chi_{\{t_{k'+1}\}}(t) \oplus ...$$

For any $i \in \{1,...,n\}$, the sets $\{t_k \mid k \geq 0, \rho_i(t_k) = 1\}, \{t_k \mid k \geq k', \rho_i(t_k) = 1\}$ are simultaneously bounded or unbounded from above.

**Remark 19** From Theorem 18 a) we get the following conclusion. For $\alpha \in \hat{\Pi}_n'$, we have the equivalence
$$\alpha \in \hat{\Pi}_n \Leftrightarrow \forall k \in \mathbf{N}, \hat{\sigma}^k(\alpha) \in \hat{\Pi}_n.$$

## 5. Flows

**Definition 20** For the function $\Phi : \mathbf{B}^n \to \mathbf{B}^n$ and $\lambda \in \mathbf{B}^n$, we define $\Phi^\lambda : \mathbf{B}^n \to \mathbf{B}^n$ by $\forall \mu \in \mathbf{B}^n$,
$$\Phi^\lambda(\mu) = (\overline{\lambda_1} \cdot \mu_1 \oplus \lambda_1 \cdot \Phi_1(\mu),...,\overline{\lambda_n} \cdot \mu_n \oplus \lambda_n \cdot \Phi_n(\mu)).$$

**Definition 21** Let $\alpha^0,..., \alpha^k, \alpha^{k+1} \in \mathbf{B}^n$, $k \geq 0$. We define the functions $\Phi^{\alpha^0...\alpha^k\alpha^{k+1}} : \mathbf{B}^n \to \mathbf{B}^n$ iteratively by $\forall \mu \in \mathbf{B}^n$,
$$\Phi^{\alpha^0...\alpha^k\alpha^{k+1}}(\mu) = \Phi^{\alpha^{k+1}}(\Phi^{\alpha^0...\alpha^k}(\mu)).$$

**Definition 22** a) The function $\mathbf{B}^n \times \mathbf{N}_- \times \hat{\Pi}_n \ni (\mu,k,\alpha) \mapsto \hat{\Phi}^\alpha(\mu,k) \in \mathbf{B}^n$ defined by $\forall k \in \mathbf{N}_-$,
$$\hat{\Phi}^\alpha(\mu,k) = \begin{cases} \mu, \text{if } k = -1, \\ \Phi^{\alpha^0...\alpha^k}(\mu), \text{if } k \geq 0 \end{cases}$$

is called (discrete time) *evolution function*, or (*state*) *transition function*, or *next state function*. $\mathbf{B}^n$ is called *state space* (or *phase space*), $\mu$ is called the *initial* (*value of the*) *state* and $\alpha$ is the *computation function*. The value
$$\hat{x}(k) = \hat{\Phi}^a(\mu,k)$$

is the state $\hat{x}(k)$ resulted at the time instant $k$ from the initial (value of the) state $\mu$ under the (action of the) computation function $\alpha$.

b) We define the function $\mathbf{B}^n \times \mathbf{R} \times \Pi_n \ni (\mu,t,\rho) \mapsto \Phi^\rho(\mu,t) \in \mathbf{B}^n$ in the following way. Let $\forall t \in \mathbf{R}$,
$$\rho(t) = \alpha^0 \cdot \chi_{\{t_0\}}(t) \oplus \alpha^1 \cdot \chi_{\{t_1\}}(t) \oplus ... \oplus \alpha^k \cdot \chi_{\{t_k\}}(t) \oplus ... \quad (12)$$

where $\alpha \in \hat{\Pi}_n$ and $(t_k) \in Seq$. Then
$$\Phi^\rho(\mu,t) = \hat{\Phi}^\alpha(\mu,-1) \cdot \chi_{(-\infty,t_0)}(t) \oplus \hat{\Phi}^\alpha(\mu,0) \cdot \chi_{[t_0,t_1)}(t) \oplus ... \oplus \hat{\Phi}^\alpha(\mu,k) \cdot \chi_{[t_k,t_{k+1})}(t) \oplus ...$$

is called (real time) *evolution function*, or (*state*) *transition function*, or *next state function*. $\mathbf{B}^n$ is the *state space*, $\mu$ is the *initial* (*value of the*) *state* and $\rho$ is the *computation function*. The value

$$x(t) = \Phi^\rho(\mu, t)$$

is the state resulted at the time instant $t$ from the initial (value of the) state $\mu$ under the (action of the) computation function $\rho$.

**Definition 23** a) We fix $\mu \in \mathbf{B}^n$ and $\alpha \in \hat{\Pi}_n$ in the argument of the discrete time evolution function. The signal $\hat{\Phi}^\alpha(\mu, \cdot) \in \hat{S}^{(n)}$ is called (discrete time) *flow* (*through* $\mu$, *under* $\alpha$) and, more general, if previously $\alpha \in \hat{\Pi}'_n$, then $\hat{\Phi}^\alpha(\mu, \cdot)$ is called *semi-flow*.

b) We fix $\mu \in \mathbf{B}^n$ and $\rho \in \Pi_n$ in the argument of the real time evolution function. The signal $\Phi^\rho(\mu, \cdot) \in S^{(n)}$ is called (real time) *flow* (*through* $\mu$, *under* $\rho$). More general, if previously $\rho \in \Pi'_n$, then $\Phi^\rho(\mu, \cdot)$ is called *semi-flow*.

**Remark 24** The function $\Phi$ applied to the argument $\mu$ is computed on all its coordinates: $\Phi(\mu) = (\Phi_1(\mu), ..., \Phi_n(\mu))$. The function $\Phi^\lambda$ applied to $\mu$ computes those coordinates $\Phi_i$ of $\Phi$ for which $\lambda_i = 1$ and it does not compute those coordinates $\Phi_i$ for which $\lambda_i = 0$: $\forall i \in \{1, ..., n\}$,

$$\Phi_i^\lambda(\mu) = \begin{cases} \Phi_i(\mu), \text{if } \lambda_i = 1, \\ \mu_i, \text{if } \lambda_i = 0. \end{cases}$$

Unlike the usual computations from the dynamical systems theory that happen synchronously on all the coordinates: $\Phi(\mu), (\Phi \circ \Phi)(\mu), (\Phi \circ \Phi \circ \Phi)(\mu), ...$ here things happen on some coordinates only, as shown in Definitions 20, 21, 22. The asynchronous flows represent a generalization of the computations from the dynamical systems theory, since the constant sequence $\alpha^k = (1, ..., 1) \in \mathbf{B}^n, k \in \mathbf{N}$ belongs to $\hat{\Pi}_n$ and it gives for any $\mu \in \mathbf{B}^n$, that $\Phi^{\alpha^0}(\mu) = \Phi(\mu)$, $\Phi^{\alpha^0 \alpha^1}(\mu) = (\Phi \circ \Phi)(\mu)$, $\Phi^{\alpha^0 \alpha^1 \alpha^2}(\mu) = (\Phi \circ \Phi \circ \Phi)(\mu), ...$

**Remark 25** We give the meaning of progressiveness: $\alpha \in \hat{\Pi}_n, \rho \in \Pi_n$ show that $\hat{\Phi}^\alpha(\mu, \cdot), \Phi^\rho(\mu, \cdot)$ compute each coordinate $\Phi_i, i = \overline{1, n}$ infinitely many times as $k \to \infty$. In electrical engineering, this corresponds to the so called *unbounded delay model of computation of the Boolean functions*, stating basically that each coordinate $i$ of $\Phi$ is computed independently on the other coordinates, in finite time.

**Remark 26** In the following we shall always suppose that the progressiveness requirement on $\alpha, \rho$ is fulfilled, thus we shall work with flows.

## 6. Consistency, composition and causality

**Remark 27** The properties stated in Theorems 28, 29, 30 and 31 to follow are the adaptation to the present context of the properties of consistency, composition and causality of the transition function that are contained in the definition of a dynamical system from [9], page 11. At the same page, the authors show that the words 'dynamical', 'non-anticipatory' and 'causal' have approximately the same meaning, making us conclude that the property of causality to be introduced may be also called non-anticipation. We must add here the remark that in the cited work the systems had an input, unlike here where it is convenient to omit this aspect, and consequently there causality referred to the input, unlike here where it refers to the computation function. The input controls the state and the computation function shows when and how the state is computed.

We suppose in this section that a function $\Phi : \mathbf{B}^n \to \mathbf{B}^n$ is given, together with $\mu \in \mathbf{B}^n, \alpha \in \hat{\Pi}_n$ and $\rho \in \Pi_n$. The relation between $\alpha$ and $\rho$ is given by (12), where $(t_k) \in Seq$.

**Theorem 28** (Consistency)

$$\hat{\Phi}(\mu, -1) = \mu, \tag{13}$$

$$\Phi^\rho(\mu, t_0 - 0) = \mu. \tag{14}$$

**Proof.**  a) This follows from Definition 22.

b) Definition 22 shows that we have

$$\forall t < t_0, \Phi^\rho(\mu,t) = \mu,$$

wherefrom (14) follows.

**Theorem 29** (Composition) a) $\forall k' \in \mathbf{N}, \forall k \in \mathbf{N}_-,$

$$\hat{\sigma}^{k'}(\hat{\Phi}^\alpha(\mu,\cdot))(k) = \hat{\Phi}^{\hat{\sigma}^{k'}(\alpha)}(\hat{\Phi}^\alpha(\mu,k'-1),k). \tag{15}$$

b) $\forall t' \in \mathbf{R}, \forall t \in \mathbf{R},$

$$\sigma^{t'}(\Phi^\rho(\mu,\cdot))(t) = \Phi^{\sigma^{t'}(\rho)}(\Phi^\rho(\mu,t'-0),t). \tag{16}$$

**Proof.** Let us notice first of all that $\alpha \in \hat{\Pi}_n, \rho \in \Pi_n \Rightarrow \hat{\sigma}^{k'}(\alpha) \in \hat{\Pi}_n, \sigma^{t'}(\rho) \in \Pi_n$ result from Theorem 18 and Remark 19, wherefrom the right members of equations (15), (16) make sense.

a) We have the following possibilities.

Case $k' = 0, k \in \mathbf{N}_-$ arbitrary, when

$$\hat{\sigma}^0(\hat{\Phi}^\alpha(\mu,\cdot))(k) = \hat{\Phi}^\alpha(\mu,k) = \hat{\Phi}^{\hat{\sigma}^0(\alpha)}(\mu,k) = \hat{\Phi}^{\hat{\sigma}^0(\alpha)}(\hat{\Phi}^\alpha(\mu,-1),k).$$

Case $k' \geq 1, k = -1$

$$\hat{\sigma}^{k'}(\hat{\Phi}^\alpha(\mu,\cdot))(-1) = \hat{\Phi}^\alpha(\mu,k'-1) = \hat{\Phi}^{\hat{\sigma}^{k'}(\alpha)}(\hat{\Phi}^\alpha(\mu,k'-1),-1).$$

Case $k' \geq 1, k \in \mathbf{N}$ arbitrary, for which

$$\hat{\sigma}^{k'}(\hat{\Phi}^\alpha(\mu,\cdot))(k) = \hat{\Phi}^\alpha(\mu,k+k') = \Phi^{\alpha^0...\alpha^{k'}\alpha^{k'+1}...\alpha^{k'+k}}(\mu)$$

$$= \Phi^{\alpha^{k'}\alpha^{k'+1}...\alpha^{k'+k}}(\Phi^{\alpha^0...\alpha^{k'-1}}(\mu)) = \hat{\Phi}^{\alpha^{k'},\alpha^{k'+1},\alpha^{k'+2},\cdots}(\Phi^{\alpha^0...\alpha^{k'-1}}(\mu),k)$$

$$= \hat{\Phi}^{\hat{\sigma}^{k'}(\alpha)}(\Phi^{\alpha^0...\alpha^{k'-1}}(\mu),k) = \hat{\Phi}^{\hat{\sigma}^{k'}(\alpha)}(\hat{\Phi}^\alpha(\mu,k'-1),k).$$

b) Equation (12) shows that we can put $\Phi^\rho(\mu,\cdot)$ under the form

$$\Phi^\rho(\mu,t) = \mu \cdot \chi_{(-\infty,t_0)}(t) \oplus \Phi^{\alpha^0}(\mu) \cdot \chi_{[t_0,t_1)}(t) \oplus ... \oplus \Phi^{\alpha^0...\alpha^k}(\mu) \cdot \chi_{[t_k,t_{k+1})}(t) \oplus ... \tag{17}$$

We take an arbitrary $t' \in \mathbf{R}$ and we have the following possibilities.

Case $t' \leq t_0$

In this situation

$$\sigma^{t'}(\rho)(t) = \rho(t),$$

$$\Phi^\rho(\mu,t'-0) = \mu,$$

thus

$$\sigma^{t'}(\Phi^\rho(\mu,\cdot))(t) = \Phi^\rho(\mu,t) = \Phi^{\sigma^{t'}(\rho)}(\Phi^\rho(\mu,t'-0),t).$$

Case $\exists k \in \mathbf{N}, t' \in (t_k, t_{k+1}]$

In this case we infer

$$\sigma^{t'}(\rho)(t) = \alpha^{k+1} \cdot \chi_{\{t_{k+1}\}}(t) \oplus \alpha^{k+2} \cdot \chi_{\{t_{k+2}\}}(t) \oplus ...$$

$$\Phi^\rho(\mu,t'-0) = \Phi^{\alpha^0...\alpha^k}(\mu),$$

$$\sigma^{t'}(\Phi^\rho(\mu,\cdot))(t) = \Phi^{\alpha^0...\alpha^k}(\mu) \cdot \chi_{(-\infty,t_{k+1})}(t) \oplus \Phi^{\alpha^0...\alpha^k\alpha^{k+1}}(\mu) \cdot \chi_{[t_{k+1},t_{k+2})}(t) \oplus ...$$

$$= \Phi^{\alpha^0...\alpha^k}(\mu) \cdot \chi_{(-\infty,t_{k+1})}(t) \oplus \Phi^{\alpha^{k+1}}(\Phi^{\alpha^0...\alpha^k}(\mu)) \cdot \chi_{[t_{k+1},t_{k+2})}(t) \oplus ...$$

$$= \Phi^{\alpha^{k+1} \cdot \chi_{\{t_{k+1}\}} \oplus \alpha^{k+2} \cdot \chi_{\{t_{k+2}\}} \oplus \cdots}(\Phi^{\alpha^0...\alpha^k}(\mu),t)$$

$$= \Phi^{\sigma^{t'}(\rho)}(\Phi^{\alpha^0...\alpha^k}(\mu),t) = \Phi^{\sigma^{t'}(\rho)}(\Phi^\rho(\mu,t'-0),t).$$

**Theorem 30** (Composition) a) For arbitrary $k' \in \mathbf{N}_-$ we can write: $\forall k \geq k'$,

$$\hat{\Phi}^\alpha(\mu,k) = \hat{\Phi}^{\hat{\sigma}^{k'+1}(\alpha)}(\hat{\Phi}^\alpha(\mu,k'),k-k'-1). \tag{18}$$

b) $\forall t' \in \mathbf{R}$ we have: $\forall t \geq t'$,

$$\Phi^\rho(\mu,t) = \Phi^{\rho \cdot \chi(t',\infty)}(\Phi^\rho(\mu,t'),t). \tag{19}$$

**Proof.** a) We make the substitution $k = k'+p$, where $p \in \mathbf{N}$ and we prove (18) by induction on $p$. For $p = 0$, (18) becomes

$$\hat{\Phi}^\alpha(\mu,k') = \hat{\Phi}^{\hat{\sigma}^{k'+1}(\alpha)}(\hat{\Phi}^\alpha(\mu,k'),-1),$$

obvious. We suppose that

$$\hat{\Phi}^\alpha(\mu,k'+p) = \hat{\Phi}^{\hat{\sigma}^{k'+1}(\alpha)}(\hat{\Phi}^\alpha(\mu,k'),p-1) \tag{20}$$

is true and we infer that

$$\hat{\Phi}^\alpha(\mu,k'+p+1) = \Phi^{\alpha^{k'+p+1}}(\hat{\Phi}^\alpha(\mu,k'+p))$$

$$\stackrel{(20)}{=} \Phi^{\alpha^{k'+p+1}}(\hat{\Phi}^{\alpha^{k'+1},\alpha^{k'+2},...\alpha^{k'+p},...}(\hat{\Phi}^\alpha(\mu,k'),p-1))$$

$$= \Phi^{\alpha^{k'+p+1}}(\Phi^{\alpha^{k'+1}\alpha^{k'+2}...\alpha^{k'+p}}(\hat{\Phi}^\alpha(\mu,k')))$$

$$= \Phi^{\alpha^{k'+1}\alpha^{k'+2}...\alpha^{k'+p}\alpha^{k'+p+1}}(\hat{\Phi}^\alpha(\mu,k'))$$

$$= \hat{\Phi}^{\alpha^{k'+1},\alpha^{k'+2},...\alpha^{k'+p},\alpha^{k'+p+1},...}(\hat{\Phi}^\alpha(\mu,k'),p) = \hat{\Phi}^{\hat{\sigma}^{k'+1}(\alpha)}(\hat{\Phi}^\alpha(\mu,k'),p).$$

b) Indeed, we shall suppose in the following that (12) is true. In the case $t' < t_0$, we have

$$\rho(t) \cdot \chi_{(t',\infty)}(t) = \rho(t),$$

$$\Phi^\rho(\mu,t') = \mu$$

and (19) is true under the form $\forall t \geq t', \Phi^\rho(\mu,t) = \Phi^\rho(\mu,t)$.

In the case $t' \in [t_k, t_{k+1}), k \in \mathbf{N}$,

$$\rho(t) \cdot \chi_{(t',\infty)}(t) = \alpha^{k+1} \cdot \chi_{\{t_{k+1}\}}(t) \oplus \alpha^{k+2} \cdot \chi_{\{t_{k+2}\}}(t) \oplus ...$$

$$\Phi^\rho(\mu,t') = \Phi^{\alpha^0...\alpha^k}(\mu),$$

$\Phi^\rho(\mu,t)$ is given by (17) and

$$\Phi^{\rho \cdot \chi(t',\infty)}(\Phi^\rho(\mu,t'),t) = \Phi^{\alpha^{k+1} \cdot \chi_{\{t_{k+1}\}} \oplus \alpha^{k+2} \cdot \chi_{\{t_{k+2}\}} \oplus ...}(\Phi^{\alpha^0...\alpha^k}(\mu),t)$$

$$= \Phi^{\alpha^0...\alpha^k}(\mu) \cdot \chi_{(-\infty,t_{k+1})}(t) \oplus \Phi^{\alpha^0...\alpha^{k+1}}(\mu) \cdot \chi_{[t_{k+1},t_{k+2})}(t) \oplus ...$$

For $t \geq t'$, (19) is true.

**Theorem 31** (Causality) For any $k \in \mathbf{N}$ and any $\alpha, \beta \in \hat{\Pi}_n$ with

$$\forall k' \in \{1,...,k\}, \alpha^{k'} = \beta^{k'}, \tag{21}$$

we have

$$\hat{\Phi}^\alpha(\mu,k) = \hat{\Phi}^\beta(\mu,k). \tag{22}$$

b) Let $t' \in \mathbf{R}$ and $\rho, \rho' \in \Pi_n$ with the property that

$$\forall t \leq t', \rho(t) = \rho'(t). \tag{23}$$

Then

$$\Phi^\rho(\mu,t') = \Phi^{\rho'}(\mu,t'). \tag{24}$$

**Proof.** a) We infer

$$\hat{\Phi}^\alpha(\mu,k) = \Phi^{\alpha^0...\alpha^k}(\mu) = \Phi^{\beta^0...\beta^k}(\mu) = \hat{\Phi}^\beta(\mu,k). \tag{25}$$

b) If $t' \in \mathbf{R}$ is such that

$$\forall t \leq t', \rho(t) = \rho'(t) = (0,...,0), \tag{26}$$

then

$$\Phi^\rho(\mu,t') = \mu = \Phi^{\rho'}(\mu,t').$$

Let $k \in \mathbf{N}$ be arbitrary and fixed. We suppose that $\alpha,\beta \in \hat{\Pi}_n$ and $(t_j),(t'_j) \in Seq$ exist such that

$$\forall k' \in \{0,...,k\}, t_{k'} = t'_{k'},$$

$$\rho(t) = \alpha^0 \cdot \chi_{\{t_0\}}(t) \oplus ... \oplus \alpha^k \cdot \chi_{\{t_k\}}(t) \oplus \alpha^{k+1} \cdot \chi_{\{t_{k+1}\}}(t) \oplus \alpha^{k+2} \cdot \chi_{\{t_{k+2}\}}(t) \oplus ...$$

$$\rho'(t) = \alpha^0 \cdot \chi_{\{t_0\}}(t) \oplus ... \oplus \alpha^k \cdot \chi_{\{t_k\}}(t) \oplus \beta^{k+1} \cdot \chi_{\{t'_{k+1}\}}(t) \oplus \beta^{k+2} \cdot \chi_{\{t'_{k+2}\}}(t) \oplus ...$$

$\alpha^0 \neq (0,...,0)$ and $t' \in [t_k,t_{k+1}) \cap [t_k,t'_{k+1})$ hold. We get

$$\Phi^\rho(\mu,t') = \Phi^{\alpha^0...\alpha^k}(\mu) = \Phi^{\rho'}(\mu,t'). \tag{27}$$

## 7. Conclusion

The flows $\hat{\Phi}^\alpha(\mu,\cdot), \Phi^\rho(\mu,\cdot)$ fulfill properties of consistency, composition and causality, as expressed by Theorems 28, 29, 30 and 31, allowing us to consider that they define dynamical systems.